\documentclass{llncs}
\usepackage{amsmath}
\usepackage{amssymb}
\usepackage{float}
\usepackage[
    backend=biber,
    style=alphabetic,
    sorting=ynt
]{biblatex}
\usepackage{hyperref}
\usepackage{cleveref}
\usepackage{tikz}
\usetikzlibrary{arrows.meta, angles, quotes, calc, intersections}
\usepackage{multirow}
\usepackage{xspace}
\usepackage{array}
\usepackage{graphicx}
\usepackage{wrapfig}
\makeatletter
\renewcommand{\fps@table}{t}
\makeatother

\makeatletter
\renewcommand\subsection{\@startsection{subsection}{2}{\z@}%
  {1.5ex \@plus .5ex \@minus .2ex}%
  {0.5ex \@plus .2ex}%
  {\normalfont\normalsize\bfseries}}
\renewcommand\paragraph{\@startsection{paragraph}{4}{\z@}%
  {1ex}%
  {-2em}%
  {\normalfont\normalsize\bfseries}}
\makeatother

\usepackage{xcolor}

\newcommand*{\sheafhom}{\mathcal{H}\kern -.5pt om}
\newcommand{\Real}{\mathbb{R}}

\newcommand{\boldalpha}{\boldsymbol{\alpha}}

\newcommand{\navBaseLOS}[1][0.7]{%
  \begin{tikzpicture}[baseline=(current bounding box.center), scale=#1]
    \node[circle, inner sep=1pt] (M)  at ( 0.0,1.0) {$M$};
    \node[circle, inner sep=1pt] (SA) at (-1.0,0.0) {$S$};
    \node[circle, inner sep=1pt] (SB) at ( 1.0,0.0) {$S$};
    \draw[dotted, line width=1pt] (M)  -- (SA);
    \draw[dotted, line width=1pt] (M)  -- (SB);
    \draw[line width=0.8pt] (SA) -- (SB);
  \end{tikzpicture}}

\newcommand{\navBaseThreeS}[1][0.7]{%
  \begin{tikzpicture}[baseline=(current bounding box.center), scale=#1]
    \node[draw, circle, inner sep=1pt] (SA) at ( 0.0,1.1) {$S$};
    \node[draw, circle, inner sep=1pt] (SB) at (-1.0,0.0) {$S$};
    \node[draw, circle, inner sep=1pt] (SC) at ( 1.0,0.0) {$S$};
    \draw[line width=0.8pt] (SA) -- (SB);
    \draw[line width=0.8pt] (SA) -- (SC);
    \draw[line width=0.8pt] (SB) -- (SC);
  \end{tikzpicture}}

\newcommand{\navBaseTwoMOneS}[1][0.7]{%
  \begin{tikzpicture}[baseline=(current bounding box.center), scale=#1]
    \node[circle, inner sep=1pt] (M1) at (-1.0,0.0) {$M$};
    \node[circle, inner sep=1pt] (M2) at ( 1.0,0.0) {$M$};
    \node[circle, inner sep=1pt] (S)  at ( 0.0,1.1) {$S$};
    \draw[line width=0.8pt] (M1) -- (S);
    \draw[line width=0.8pt] (M2) -- (S);
    \draw[line width=0.8pt] (M1) -- (M2);
  \end{tikzpicture}}

\newcommand{\navBaseMTwoS}[1][0.7]{%
  \begin{tikzpicture}[baseline=(current bounding box.center), scale=#1]
    \node[circle, inner sep=1pt] (M)  at ( 0.0,1.0) {$M$};
    \node[circle, inner sep=1pt] (SA) at (-1.0,0.0) {$S$};
    \node[circle, inner sep=1pt] (SB) at ( 1.0,0.0) {$S$};
    \draw[line width=0.8pt] (M)  -- (SA);
    \draw[line width=0.8pt] (M)  -- (SB);
    \draw[line width=0.8pt] (SA) -- (SB);
  \end{tikzpicture}}

\newcommand{\navBaseMTwoSSame}[1][0.7]{%
  \begin{tikzpicture}[baseline=(current bounding box.center), scale=#1]
    \node[circle, inner sep=1pt] (M)  at ( 0.0,1.0) {$M$};
    \node[draw, dotted, circle, inner sep=1pt] (SA) at (-1.0,0.0) {$S$};
    \node[draw, dotted, circle, inner sep=1pt] (SB) at ( 1.0,0.0) {$S$};
    \draw[line width=0.8pt] (M)  -- (SA);
    \draw[line width=0.8pt] (M)  -- (SB);
    \draw[line width=0.8pt] (SA) -- (SB);
  \end{tikzpicture}}

\newcommand{\navBaseMTwoSRangeRate}[1][0.7]{%
  \begin{tikzpicture}[baseline=(current bounding box.center), scale=#1]
    \node[circle, inner sep=1pt] (M)  at ( 0.0,1.1) {$M$};
    \node[circle, inner sep=1pt] (SA) at (-1.1,0.0) {$S$};
    \node[circle, inner sep=1pt] (SB) at ( 1.1,0.0) {$S$};
    \draw[double, double distance=2pt, line width=0.4pt] (M)  -- (SA);
    \draw[double, double distance=2pt, line width=0.4pt] (M)  -- (SB);
    \draw[double, double distance=2pt, line width=0.4pt] (SA) -- (SB);
  \end{tikzpicture}}

\newcommand{\navBaseSixS}[1][0.7]{%
  \begin{tikzpicture}[baseline=(current bounding box.center), scale=#1]
    \node[circle, inner sep=1pt] (A) at ( 0.00, 2.0) {$S$};
    \node[circle, inner sep=1pt] (B) at (-1.35, 1.0) {$S$};
    \node[circle, inner sep=1pt] (C) at ( 1.35, 1.0) {$S$};
    \node[circle, inner sep=1pt] (D) at ( 0.00,-1.2) {$S$};
    \node[circle, inner sep=1pt] (E) at (-1.50,-0.9) {$S$};
    \node[circle, inner sep=1pt] (F) at ( 1.50,-0.9) {$S$};
    \draw[line width=0.8pt] (A) -- (B);
    \draw[line width=0.8pt] (A) -- (C);
    \draw[line width=0.8pt] (A) -- (D);
    \draw[line width=0.8pt] (B) -- (C);
    \draw[line width=0.8pt] (B) -- (D);
    \draw[line width=0.8pt] (C) -- (D);
    \draw[line width=0.8pt] (A) -- (E);
    \draw[line width=0.8pt] (B) -- (E);
    \draw[line width=0.8pt] (C) -- (E);
    \draw[line width=0.8pt] (A) -- (F);
    \draw[line width=0.8pt] (B) -- (F);
    \draw[line width=0.8pt] (C) -- (F);
  \end{tikzpicture}}

\newcommand{\navBaseThreeSTopKnownCbottomSSame}[1][0.7]{%
  \begin{tikzpicture}[baseline=(current bounding box.center), scale=#1]
    \node[draw, circle, inner sep=1pt] (SA) at ( 0.0, 1.1) {$S$};
    \node[draw, dotted, circle, inner sep=1pt] (SB) at (-1.0,0.0) {$S$};
    \node[draw, dotted, circle, inner sep=1pt] (SC) at ( 1.0,0.0) {$S$};
    \draw[line width=0.8pt, transform canvas={yshift=1.5pt}] (SA) -- (SB);
    \draw[dotted, line width=0.8pt, transform canvas={yshift=-1.5pt}] (SA) -- (SB);
    \draw[line width=0.8pt, transform canvas={yshift=1.5pt}] (SA) -- (SC);
    \draw[dotted, line width=0.8pt, transform canvas={yshift=-1.5pt}] (SA) -- (SC);
  \end{tikzpicture}}

\newcommand{\navBaseThreeSAllUnknownC}[1][0.7]{%
  \begin{tikzpicture}[baseline=(current bounding box.center), scale=#1]
    \node[circle, inner sep=1pt] (SA) at ( 0.0, 1.1) {$S$};
    \node[circle, inner sep=1pt] (SB) at (-1.0,0.0) {$S$};
    \node[circle, inner sep=1pt] (SC) at ( 1.0,0.0) {$S$};
    \draw[line width=0.8pt, transform canvas={yshift=1.5pt}] (SA) -- (SB);
    \draw[dotted, line width=0.8pt, transform canvas={yshift=-1.5pt}] (SA) -- (SB);
    \draw[line width=0.8pt, transform canvas={yshift=1.5pt}] (SA) -- (SC);
    \draw[dotted, line width=0.8pt, transform canvas={yshift=-1.5pt}] (SA) -- (SC);
  \end{tikzpicture}}

\newcommand{\navBaseTwoSOneCRangeLOS}[1][0.7]{%
  \begin{tikzpicture}[baseline=(current bounding box.center), scale=#1]
    \node[draw, circle, inner sep=1pt] (SA) at (-1.0, 0.0) {$S$};
    \node[circle, inner sep=1pt] (SB) at ( 1.0, 0.0) {$S$};
    \draw[line width=0.8pt, transform canvas={yshift=1.5pt}] (SA) -- (SB);
    \draw[dotted, line width=0.8pt, transform canvas={yshift=-1.5pt}] (SA) -- (SB);
  \end{tikzpicture}}

\newcommand{\navBaseTwoSRangeLOS}[1][0.7]{%
  \begin{tikzpicture}[baseline=(current bounding box.center), scale=#1]
    \node[draw, circle, inner sep=1pt] (SA) at (-1.0, 0.0) {$S$};
    \node[draw, circle, inner sep=1pt] (SB) at ( 1.0, 0.0) {$S$};
    \draw[line width=0.8pt, transform canvas={yshift=1.5pt}] (SA) -- (SB);
    \draw[dotted, line width=0.8pt, transform canvas={yshift=-1.5pt}] (SA) -- (SB);
  \end{tikzpicture}}

\newcommand{\mpGraphLOSTwice}[1][0.7]{$\left(\vcenter{\hbox{\navBaseLOS[#1]}}\right)^{2}$}
\newcommand{\mpGraphThreeS}[1][0.7]{\navBaseThreeS[#1]}
\newcommand{\mpGraphTwoMOneS}[1][0.7]{\navBaseTwoMOneS[#1]}
\newcommand{\mpGraphMTwoSSame}[1][0.7]{\navBaseMTwoSSame[#1]}
\newcommand{\mpGraphMTwoSDistinctTwice}[1][0.7]{$\left(\vcenter{\hbox{\navBaseMTwoS[#1]}}\right)^{2}$}
\newcommand{\mpGraphMTwoSRangeRate}[1][0.7]{\navBaseMTwoSRangeRate[#1]}
\newcommand{\mpGraphSixS}[1][0.7]{\navBaseSixS[#1]}
\newcommand{\mpGraphTwoSRangeLOS}[1][0.7]{\navBaseTwoSRangeLOS[#1]}
\newcommand{\mpGraphTwoSOneCRangeLOS}[1][0.7]{\navBaseTwoSOneCRangeLOS[#1]}
\newcommand{\mpGraphThreeSTopKnownCBottomSSame}[1][0.7]{\navBaseThreeSTopKnownCbottomSSame[#1]}
\newcommand{\mpGraphThreeSAllUnknownC}[1][0.7]{\navBaseThreeSAllUnknownC[#1]}

\title{Approximating Periodic Orbits with Algebraic Curves and Related Minimal Problems}

\author{Ruiqi Huang \and Anton Leykin\thanks{The work was partially supported by AFOSR award FA95502310512 and NSF DMS award 2001267.}}
\institute{Georgia Institute of Technology, Atlanta GA 30332, USA}

\begin{document}

\maketitle

\begin{abstract}
The Circular Restricted Three-Body Problem (CR3BP) models the motion of a massless
body under the gravitational influence of two primaries. We present a method for
approximating a given family of periodic orbits by low-degree implicit algebraic curves,
producing one-parameter families of algebraic orbit models.

These models enable the construction of minimal problems motivated by liaison navigation,
where spacecraft states are inferred from inter-spacecraft measurements. Relevant
applications include initial orbit determination and spacecraft positioning.

Each minimal problem defines a parameterized family of instances; for generic parameters,
the number of solutions equals the degree of the associated branched covering map. We
compute these degrees using both symbolic and numerical methods, and we outline a
homotopy-continuation-based solver construction that can be practical for low-degree cases.
\end{abstract}

\section{Introduction}

The Circular Restricted Three-Body Problem (CR3BP) seeks to determine the motion of a
massless body subject to the gravitational attraction of two massive primaries (e.g.,
Earth and Moon) \cite{szebehely_theory_1967}. The general three-body problem admits no
closed-form solution \cite{poincare_methodes_1892}; however, in the restricted setting, 
one finds rich families of periodic orbits, including the planar Lyapunov orbits and the
three-dimensional Halo orbits around the collinear Lagrange points.

In this article, for a known family of periodic orbits,
we approximate data points on a given orbit with an algebraic curve of low degree. We then fit the data from multiple orbits in the family, obtaining a model depending on one parameter linked to the orbit's energy.
The main contribution of this work is to use the above models to formulate \emph{minimal problems} for liaison navigation that use crosslink range, range rate, and line-of-sight measurements. For each minimal problem, we determine the \emph{algebraic degree}, which is equal to the number of complex solutions in the generic case. 

The problems with small algebraic degrees admit efficient solvers using parameter homotopy continuation, which could quickly obtain approximate solutions that can be used as initial guesses for dynamics-aware methods such as the liaison navigation approach in~\cite{liaison:keric-born}.

\section{Preliminaries}

\subsection{Rotating Frame and CR3BP}
We work in the \emph{rotating (synodic) frame} in which the two primaries of masses $1-\mu$
and $\mu$ are fixed at $(-\mu,0)$ and $(1-\mu,0)$ on the $x$-axis, and the frame rotates around
the $z$-axis at unit angular velocity. 

{

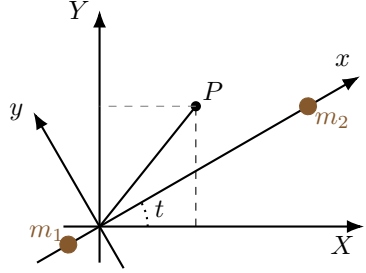
\begin{wrapfigure}{l}{0.45\linewidth}
    \centering
    \begin{tikzpicture}[scale=1.6, >=Latex]
\draw[->, thick] (-0.3,0) -- (2.2,0) node[below left] {$X$};
\draw[->, thick] (0,-0.3) -- (0,1.8) node[left] {$Y$};

\begin{scope}[rotate=30]
  \draw[->, thick] (-0.6,0) -- (2.5,0) node[above left] {$x$};
  \draw[->, thick] (0,-0.4) -- (0,1.1) node[left] {$y$};
\end{scope}

\filldraw[brown!70!black, rotate around={30:(0,0)}] (-0.3,0) circle (2pt) node[above left=-2pt] {$m_1$};
\filldraw[brown!70!black, rotate around={30:(0,0)}] (2,0) circle (2pt) node[below right=-1pt] {$m_2$};

\coordinate (P) at (0.8,1.0);
\fill (P) circle (1.2pt) node[above right=-1pt] {$P$};

\draw[dashed] (0.8,0) |- (P);
\draw[thick] (0,0) -- (P);

\draw[thick, dotted] (0.4,0) arc[start angle=0, end angle=30, radius=0.4];
\node at (0.50,0.15) {$t$};

\draw[dashed, gray] (P) -- (0,1);
    \end{tikzpicture}
    \caption{Inertial $(X,Y)$ and rotating frame $(x,y)$ in the planar case.}
    \label{fig:rotating_frame}
\end{wrapfigure}
\noindent The state $(x,y,z,\dot x,\dot y,\dot z)$ of the massless body is
guided by the equations of motion:
\begin{align} \label{eq:ODE}
    \ddot{x} &= 2\dot{y}+x \nonumber\\
             &\quad -\frac{(1-\mu)(x+\mu)}{r_1^3}-\frac{\mu(x-(1-\mu))}{r_2^3},\\
    \ddot{y} &= -2\dot{x}+y-\frac{(1-\mu)y}{r_1^3}-\frac{\mu y}{r_2^3}, \nonumber\\
    \ddot{z} &= -\frac{(1-\mu)z}{r_1^3}-\frac{\mu z}{r_2^3}, \nonumber
\end{align}
where $r_1^2=(x+\mu)^2+y^2+z^2$ and $r_2^2=(x-(1-\mu))^2+y^2+z^2$.

}
\smallskip
We imagine the rotating frame sharing the center of mass of the two primaries with the inertial frame where angle $t$ in \Cref{fig:rotating_frame} changes linearly with time, so the dynamics \eqref{eq:ODE} involve Coriolis and centrifugal forces. 

\subsection{Jacobi Constant and Periodic Orbits}

The system \eqref{eq:ODE} is Hamiltonian; the conserved quantity is the \emph{Jacobi
constant} $C = -2H_r$, expressed explicitly as
\begin{equation} \label{eq:Jacobi}
    C = (x^2+y^2) + \frac{2(1-\mu)}{r_1} + \frac{2\mu}{r_2} - (\dot{x}^2+\dot{y}^2)
      = 2\Omega - (\dot{x}^2+\dot{y}^2),
\end{equation}
where the \emph{effective potential} is $\Omega = \tfrac{1}{2}(x^2+y^2)+\tfrac{1-\mu}{r_1}+\tfrac{\mu}{r_2}$.
Each family of periodic orbits that we consider is parameterized by $C$.

\emph{Lyapunov orbits} are planar periodic solutions around the Lagrange points
$L_1$ or $L_2$.
\emph{Halo orbits} are three-dimensional periodic solutions with an out-of-plane component;
they also exist around $L_1$ and $L_2$. Both families are continuous in $C$. In our numerical experiments, we use orbit
data from the JPL Three-Body Periodic Orbit Catalog \cite{nasa_three-body_nodate} for
the Earth--Moon system ($\mu\approx 0.01215$).

\section{Fitting Algebraic Curves}

\subsection{Polynomial Model for Lyapunov and Halo orbits}

Each \emph{Lyapunov} orbit is symmetric about the $x$-axis; only even powers of $y$ appear in
an implicit polynomial equation that we use as our \emph{quartic model}:
\begin{equation} \label{eq:Lyapunov-model}
    g(x,y) := \alpha_1 x + \alpha_2 x^2 + \alpha_3 x^3 + \alpha_4 x^4
             + \alpha_5 y^2 + \alpha_6 xy^2 + \alpha_7 x^2y^2 + \alpha_8 y^4 = 1,
\end{equation}
with 8 free coefficients $\boldalpha=(\alpha_1,\dots,\alpha_8)$.

For \emph{Halo} orbits, change the coordinate frame to the basis
$\mathbf{u}=(-\tfrac{1}{\sqrt 2},0,\tfrac{1}{\sqrt 2})$,
$\mathbf{v}=(0,1,0)$,
$\mathbf{w}=(\tfrac{1}{\sqrt 2},0,\tfrac{1}{\sqrt 2})$
and fit \emph{two} equations in the $(u,v,w)$-coordinates:
\begin{equation} \label{eq:Halo-model}
    g(u,v)=\alpha_1 u + \alpha_2 u^2 + \cdots + \alpha_8 v^4 = 1, \qquad
    h(u,v)=\beta_1 + \beta_2 u + \cdots + \beta_9 v^4 = w,
\end{equation}
effectively reducing the model to a planar curve: the first equation describes the projected curve in the $uv$-plane, and the second
gives the third coordinate $w$ as a polynomial in $(u,v)$.
\subsection{Fitting One Orbit}

\begin{wrapfigure}{r}{0.25\textwidth}
  \centering
  \begin{tikzpicture}
    \node[anchor=south west, inner sep=0] (img) at (0,0)
      {\includegraphics[width=\linewidth]{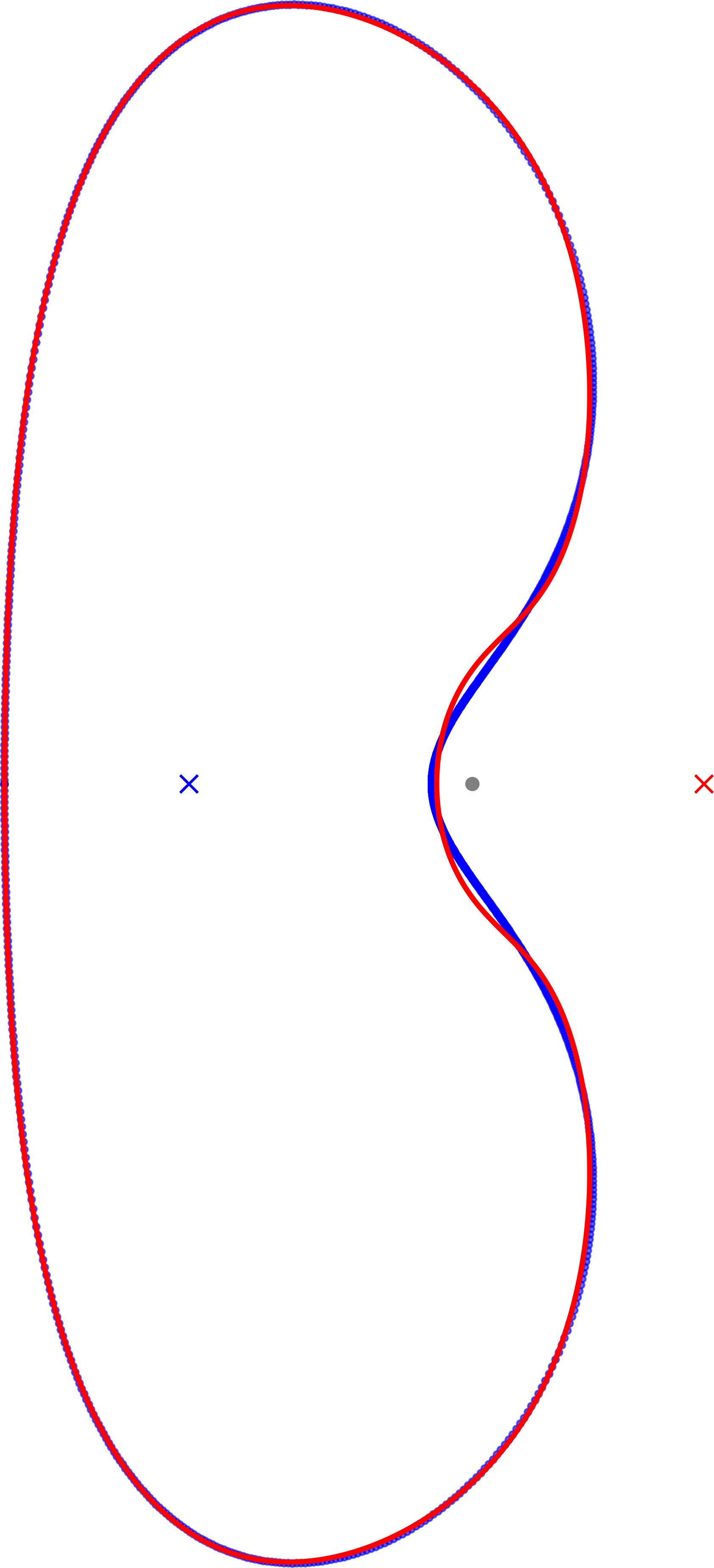}};
    \begin{scope}[x={(img.south east)},y={(img.north west)}]
      \node at (0.33,0.47) {\footnotesize $L_1$};
      \node at (0.79,0.475) {\footnotesize Moon};
      \node at (1.05,0.47) {\footnotesize $L_2$};
    \end{scope}
  \end{tikzpicture}
\end{wrapfigure}

Given data points $\{(x_i,y_i)\}_{i=1}^N$ on a single orbit (\emph{fixed} $C$), the
coefficients $\boldalpha$ are determined by the least-squares normal equations
\[
    (A^\top A)\boldalpha = A^\top \mathbf{1},
\]
where $A_{ij}=\phi_j(x_i,y_i)$ and $\phi_j$ are the monomials of~$g$.

The \emph{root mean square error} has a unique minimum at the least-squares solution. 
However, a more geometrically meaningful metric is the \emph{mean Euclidean distance}
from the data points to the fitted curve, which takes more computational effort to approximate, let alone optimize.  

We note that optimizing the (approximate) mean geometric distance is a nonlinear problem, which can be solved by a local method (e.g., gradient descent) using the initial approximation from least-squares. 

\subsection{Fitting a Family of Orbits}

The simplest tool to use for fitting a \emph{model} $g(x,y,C)\in\Real[C][x,y]$ for Lyapunov or $g,h\in\Real[C][u,v]$ for Halo is, again, the least-squares method. 

\paragraph{Two-stage procedure.}
For each orbit in a training set $\{C_k\}$, solve the least-squares problem to obtain
coefficients $\boldalpha(C_k)$. Then fit each component $\alpha_j(C)$ independently as a
low-degree polynomial in $C$ (we use cubics in experiments). For practical purposes, we partition $[C_{\min},C_{\max}]$ into subintervals
that contain data for about the same number of orbits in the JPL catalog.

\paragraph{One-stage procedure.}
Treat all orbit data simultaneously with a single least-squares solve: for instance, in the case of Lyapunov orbits, expand $g(x,y,C)$
in monomials $\{\phi_j(x,y,C)\}$ jointly in all three variables, assemble the design matrix
$A_{ij}=\phi_j(x_i,y_i,C_i)$, and solve $(A^\top A)\boldalpha = A^\top\mathbf{1}$
for the coefficient vector $\boldalpha$.

In our numerical experiments with the Lyapunov family, the one-stage procedure doesn't provide significant improvement over the two-stage for the quartic model. However, the performance of the two-stage for the sextic model is significantly poorer than the one-stage.

\section{Minimal Problems for Navigation}

A \emph{minimal problem} is a square system of parametric polynomial equations (same number of
equations and unknowns) that has finitely many complex solutions for generic values of
the parameters. The \emph{degree} of a minimal problem is the number of complex solutions
in the generic case, which equals the degree of the corresponding branched covering map (from the so-called solution variety to the parameter space).

\subsection{Descriptive notation and minimal problem examples}
In \emph{liaison navigation} \cite{liaison:keric-born}, spacecraft exchange
\emph{crosslink range} (distance) measurements $d_{AB}$ and possibly
\emph{range-rate} $\dot d_{AB}$ between point $A$ and $B$; a mothership (e.g., a moon-based or earth-based) observer may also supply
\emph{line-of-sight} (LOS) directions. We denote a \emph{mothership}, a body with a known state, by $M$ and a \emph{spacecraft}, a body with an unknown state, by
$S$; edges in the measurement graph indicate known measurements, with edge type encoding
the measurement type (see \Cref{tab:Lyapunov-Degrees}).


To give an example, we describe in some detail the formulation of two minimal problems  involving one mothership and two spacecraft (M2S). 

\paragraph{M2S (Same Orbit).}
\label{sec:example-M1S2}
\begin{wrapfigure}{r}{0pt}
  \mpGraphMTwoSSame \end{wrapfigure}
Consider a single mothership $M$ at a known position $(x_M,y_M)$, and two spacecraft $A$ and $B$
on the \emph{same} Lyapunov orbit with an unknown Jacobi constant $C$. The available
measurements are the mothership–spacecraft distances $d_{AM}$, $d_{BM}$, and the
inter-spacecraft distance $d_{AB}$. 

The unknowns are $(C, x_A, y_A, x_B, y_B)$ — five in total. The five equations are:
\begin{alignat}{2}
    &\text{(orbit model):} & \quad g(x_A,y_A,C) &= 1,\label{eq:ex1}\\
    &                      & g(x_B,y_B,C) &= 1,\label{eq:ex2}\\
    &\text{(M-crosslink):} &
      (x_A-x_M)^2+(y_A-y_M)^2 &= d_{AM}^2,\label{eq:ex3}\\
    & & (x_B-x_M)^2+(y_B-y_M)^2 &= d_{BM}^2,\label{eq:ex4}\\
    &\text{(crosslink):} &
      (x_A-x_B)^2+(y_A-y_B)^2 &= d_{AB}^2.\label{eq:ex5}
\end{alignat}
This system defines a minimal problem. The degree is \textbf{84} for the quartic model
\eqref{eq:Lyapunov-model} and \textbf{132} for the sextic analog. Computed degrees\footnote{Degrees are computed using the Gr\"{o}obner basis engine of \texttt{Macaulay2}~\cite{M2}
applied to generic problem instances; for practical
efficiency, computations are performed over a large finite field.
See \url{https://github.com/Rick3yHuang/Navigation-Problems-in-CR3BP}.
} 
for this and other minimal problems in the article are collected in \Cref{tab:Lyapunov-Degrees}. 

\paragraph{M2S (Range and Range Rate).}
\label{sec:example-rangeRate}
\begin{wrapfigure}{r}{0pt}
  \mpGraphMTwoSRangeRate \end{wrapfigure}
Consider a mothership $M$ at a known position and state, and two spacecraft $A$ and $B$ on
\emph{distinct} Lyapunov orbits with unknown Jacobi constants $C_A$ and $C_B$. Each edge in
the measurement graph carries both the \emph{range} (distance) and \emph{range rate}
(the first derivative of distance).

The unknowns are the positions and velocities of both spacecraft, together with their Jacobi
constants: $(C_A, x_A, y_A, \dot x_A, \dot y_A, C_B, x_B, y_B, \dot x_B, \dot y_B).$ The ten equations comprise three range constraints, three range-rate
constraints, and two orbit-model constraints (one per spacecraft), with the velocities
further related to $C$ via \eqref{eq:Jacobi}. The symbolic computation exceeded the time
limit (TLE) for both the quartic and sextic models.

\subsection{Minimal Problems for Halo Orbits}

Halo orbits are three-dimensional, allowing for more independent distance measurements due to less rigidity. Some Lyapunov minimal problems become non-minimal for Halo orbits, and vice versa, allowing new Halo minimal problems.
The degrees of the problems considered above that remain Halo-minimal appear in Table~\ref{tab:Lyapunov-Degrees}.

\paragraph{6S (an independent set of ranges).}
\begin{wrapfigure}{r}{0pt}
  \mpGraphSixS \end{wrapfigure}
The six-spacecraft problem is new: each spacecraft lies on a Halo orbit
with an unknown Jacobi constant, contributing two orbit-model equations. With a suitable
crosslink measurements
graph on the six nodes,\footnote{Due to rigidity in 3D space, the distances corresponding to the edges absent in the graph of the problem are algebraically dependent on the set of distances corresponding to the edges that are present.}  
the system becomes square
and defines a minimal problem. By contrast, with only the range measurement, for Lyapunov orbits, no number of spacecraft or
time instances suffices to form a minimal problem without a mothership, since each
spacecraft contributes only one orbit-model constraint and the system is always
underdetermined.

\subsection{Two spacecraft scenario and line-of-sight measurements}
We note that, in theory, the problem with sufficiently many crosslink measurements
identifies a generic trajectory in CR3BP, as
exploited in~\cite{liaison:keric-born}. In our setting, one can formulate a minimal problem involving sufficiently many higher derivatives of the range function (which could be numerically approximated from range and range rate measurements). E.g., for Lyapunov models, we need range and its five derivatives: indeed, for 10 unknowns---$C$, $x$, $y$, $\dot x$, $\dot y$ for each spacecraft---we then have 10 constraints (2 orbit constraints + 2 Jacobi constant constraints + 6 range-function constraints for orders from 0 to 5). However, such a problem has a very high algebraic degree, thus limiting its practical utility.

Another type of measurement we may introduce is the \emph{line-of-sight (LOS) direction} from one body to another; this is practically achieved through angular measurements with the distant (known) stars. LOS measurements and their applications have been studied, for instance, in~\cite{mancini2025geometric} in the angles-only setting for the two-body problem. 

LOS measurements are used either to derive additional constraints or to reduce the number of unknowns: if the direction vector from one body to another is known, the position of the second body is determined by the range.

\paragraph{2S (known range and LOS).}
\label{sec:example-2S-rangeLOS}
\begin{wrapfigure}{r}{0pt}
  \mpGraphTwoSRangeLOS \end{wrapfigure}
Consider two spacecraft $A$ and $B$ on Lyapunov orbits with \emph{known} Jacobi constants
$C_A$ and $C_B$. The available measurements are the inter-spacecraft range $d_{AB}$ and the
line-of-sight direction $\theta_{AB}$ from $A$ to $B$. Since both measurements are known,
$B$'s position is determined by $A$'s:
\[
  \begin{pmatrix} x_B \\ y_B \end{pmatrix}
  =
  \begin{pmatrix} x_A \\ y_A \end{pmatrix}
  -
  d_{AB}\begin{pmatrix} \cos\theta_{AB} \\ \sin\theta_{AB} \end{pmatrix}.
\]
The only unknowns are $(x_A, y_A)$, and the two orbit-model constraints
\begin{align}
  g(x_A,\, y_A,\, C_A) &= 1, \nonumber\\
  g(x_B,\, y_B,\, C_B) &= 1 \nonumber
\end{align}
form a minimal problem. The degree of this problem is  \textbf{16} for the quartic model and \textbf{36} for the sextic model.

Note that, while many of the minimal problems we formulated for Lyapunov work for Halo models as well, this doesn't hold here. Since we have a 3D direction vector, it introduces a constraint on $z$-coordinates, which are already determined by the second equation in \eqref{eq:Halo-model}, making the system overdetermined. 
\begin{wrapfigure}{r}{0pt}
  \mpGraphTwoSOneCRangeLOS \end{wrapfigure}
Relaxing the assumption on Jacobi constants---one known and one unknown---
we form a Halo minimal problem of degrees \textbf{96} for the quartic and \textbf{216} for the sextic model.

\paragraph{3S (known range and LOS from one spacecraft).}
\begin{wrapfigure}{r}{0pt}
  \mpGraphThreeSTopKnownCBottomSSame
\end{wrapfigure}
Assume spacecraft $A$, $B$, and $D$ lie on two Lyapunov orbits. $A$ lies on an orbit with a \emph{known} Jacobi constant $C_A$, while $B$ and $D$ lie on the \emph{same} orbit with an \emph{unknown} Jacobi constant $C_{BD}$. The available measurements are the inter-spacecraft ranges $d_{AB}$ and $d_{AD}$, together with the line-of-sight directions ${\theta}_{AB}$ and ${\theta}_{AD}$. As in the previous example, the positions of $B$ and $D$ can then be recovered from the position of $A$.
The unknowns are $(x_A,y_A)$ and $C_{BD}$, and the constraints are the three orbit-model constraints. This minimal problem has a degree of \textbf{84} for the quartic model and \textbf{198} for the sextic model.

\begin{wrapfigure}{r}{0pt}
  \mpGraphThreeSAllUnknownC
\end{wrapfigure}
As in the previous example, the Lyapunov-minimal problem must be relaxed to become Halo-minimal: assuming all three Jacobi constants are unknown, the degree is \textbf{3024} for the quartic model.

\paragraph{M2S (LOS, Two Time Instances).}
\label{sec:example-LOS2}
\begin{wrapfigure}{r}{0pt}
  \mpGraphLOSTwice \end{wrapfigure}
Consider a mothership $M$ at a known position and two spacecraft $A$, $B$ on \emph{distinct}
Lyapunov orbits with unknown Jacobi constants $C_A$, $C_B$. At each of the two time instances
$i=1,2$, the mothership measures the line-of-sight directions to both spacecraft and the
inter-spacecraft crosslink distance $d_{AB}^i$.
Since the LOS directions from $M$ are known, each spacecraft position at time $i$ is determined by its distance from $M$, leaving six unknowns:
$(C_A,\, C_B,\, d_{MA}^1,\, d_{MA}^2,\, d_{MB}^1,\, d_{MB}^2)$.
The six equations consist of four orbit-model constraints (one per spacecraft per time instance) and two crosslink constraints (one per time instance).

\section{Discussion and Future Work}

\paragraph{Other families of orbits.}
Our fitting method extends beyond Lyapunov and Halo orbits to any periodic or quasi-periodic CR3BP family and, more broadly, to other dynamical models. Being agnostic to the underlying dynamics, it suits scenarios where a \emph{timeless} (angles-only) description is preferred.

\paragraph{Homotopy continuation solvers.}
For low-degree problems (e.g., degrees 6, 84, 132), efficient solvers can be built in the spirit of numerical nonlinear algebra~\cite{bates2023numerical}: pre-compute a \emph{start system} for a generic instance, then track solution paths via parameter homotopy to new instances.

\paragraph{Subdivision and low-degree approximation.}
The degree growth observed for harder problems suggests using piecewise approximations (e.g., trigonometric polynomials in polar coordinates, or splines) to keep algebraic complexity manageable.

\paragraph{Nonlinear optimization.}
Solutions to the minimal problems serve as initial guesses for local nonlinear optimizers (e.g., Levenberg--Marquardt) that minimize equation residuals.

\paragraph{Initialization for dynamics-aware formulations.}
The algebraic orbit model is an approximation; a full navigation solution must account for the true dynamics \eqref{eq:ODE}. Minimal-problem solutions provide high-quality initial guesses for statistical estimators such as the Kalman filter~\cite{liaison:keric-born}, which refine the orbit estimate using subsequent measurements.

\printbibliography

\bigskip

\begin{table}[h]
\centering

\begin{minipage}[c]{0.68\linewidth}
\centering
\small
\begin{tabular}{|c|*{4}{>{\centering\arraybackslash}p{1.15cm}|}}
\hline
\multirow{2}{*}{\textbf{Problems}}
& \multicolumn{2}{c|}{\textbf{Lyapunov}}
& \multicolumn{2}{c|}{\textbf{Halo}} \\ \cline{2-5}
& quartic & sextic & quartic & sextic \\ \hline

\rule{0pt}{20pt}\mpGraphTwoMOneS\rule[-20pt]{0pt}{20pt}
& 6 & 6 & 48 & 72
\\ \cline{1-5}

\rule{0pt}{20pt}\mpGraphTwoSRangeLOS\rule[-20pt]{0pt}{20pt}
& 16 & 36 & --- & ---
\\ \cline{1-5}

\rule{0pt}{20pt}\mpGraphTwoSOneCRangeLOS\rule[-20pt]{0pt}{20pt}
& --- & --- & 96 & 216
\\ \cline{1-5}

\rule{0pt}{20pt}\mpGraphMTwoSSame\rule[-20pt]{0pt}{20pt}
& 84 & 132 & 2976 & \textbf{TLE}
\\ \cline{1-5}

\rule{0pt}{20pt}\mpGraphThreeSTopKnownCBottomSSame\rule[-20pt]{0pt}{20pt}
& 84 & 198 & --- & ---
\\ \cline{1-5}

\rule{0pt}{20pt}\mpGraphThreeSAllUnknownC\rule[-20pt]{0pt}{20pt}
& --- & --- & 3024 & \textbf{TLE}
\\ \cline{1-5}

\rule{0pt}{20pt}\mpGraphThreeS\rule[-20pt]{0pt}{20pt}
& 256 & 864 & 8192 & \textbf{TLE}
\\ \cline{1-5}

\rule{0pt}{26pt}\mpGraphLOSTwice\rule[-20pt]{0pt}{20pt}
& 1152 & 2592 & --- & ---
\\ \hline
\end{tabular}
\end{minipage}%
\hfill
\begin{minipage}[c]{0.28\linewidth}
\small
\begin{tabular}{@{}cl@{}}
  \tikz[baseline=-0.5ex] \draw[line width=0.8pt] (0,0) -- (0.9,0); & range \\[2pt]
  \tikz[baseline=-0.5ex] \draw[double, double distance=2pt, line width=0.4pt] (0,0) -- (0.9,0); & range and range rate \\[2pt]
  \tikz[baseline=-0.5ex] \draw[dotted, line width=1pt] (0,0) -- (0.9,0); & line-of-sight \\[2pt]
  \tikz[baseline=-0.5ex] { \draw[line width=0.8pt, transform canvas={yshift=1.5pt}] (0,0) -- (0.9,0); \draw[dotted, line width=0.8pt, transform canvas={yshift=-1.5pt}] (0,0) -- (0.9,0); } & range and LOS \\[2pt]
  \tikz[baseline=(n.base)] \node[draw, circle, inner sep=1pt] (n) {$S$}; & known $C$ \\[2pt]
  \tikz[baseline=(n.base)] \node[draw, dotted, circle, inner sep=1pt] (n) {$S$}; & same orbit \\[8pt]
  \textbf{TLE} & Time Limit \\
               & Exceeded \\
               & (${>}5$ hours) \\[8pt]
  \multicolumn{2}{@{}c@{}}{%
    \fbox{\begin{minipage}{0.95\linewidth}\centering
      \textit{All models give TLE:}\\[4pt]
      \mpGraphSixS[0.5]\\[4pt]
      \mpGraphMTwoSDistinctTwice[0.5]\\[4pt]
      \mpGraphMTwoSRangeRate[0.5]
    \end{minipage}}
  } \\
\end{tabular}
\end{minipage}

\vspace{0.8em}
\caption{Degrees of the minimal problems. Local dimension test (Jacobian has full rank at a generic point) is passed by every problem with TLE confirming that these are minimal.}

\label{tab:Lyapunov-Degrees}
\end{table}

\end{document}